\documentclass[11pt]{article}

 \usepackage{setspace}
\usepackage{natbib}

\usepackage{latexsym}
\usepackage{calc}

\usepackage{amssymb}
\usepackage{amsmath}
\usepackage{graphicx}

\usepackage{psfrag}

\usepackage{float}

\usepackage{color}

\usepackage{nicefrac}

\usepackage{ntheorem}

\usepackage{hyperref}

\usepackage{tikz}
\usetikzlibrary{backgrounds}

\newcounter{hours}\newcounter{minutes}

\def\nr{\par }

\def\beq{\begin{equation}}
\def\eeq{\end{equation}}

\newcommand{\N}{\mathbb{N}}

\newtheorem{theorem}{Theorem}

\newtheorem{lemma}{Lemma}
\newtheorem{corollary}{Corollary}

\newtheorem{proposition}{Proposition}
\newtheorem{assumption}{Assumption}
\newtheorem{definition}{Definition}

\newtheorem{example}{Example}
\newtheorem{remark}{Remark}
\newcommand{\proof}{\bf Proof: \rm \nr}

\newcommand{\qed}{\hfill $\Box$ \nr \medskip}

\def\ba{\begin{array}}
\def\ea{\end{array}}
\def\beann{\begin{eqnarray*}}
\def\eeann{\end{eqnarray*}}
\def\bea{\begin{eqnarray}}
\def\eea{\end{eqnarray}}

\def\BT{\begin{theorem}}
\def\ET{\end{theorem}}
\def\BL{\begin{lemma}}
\def\EL{\end{lemma}}
\def\BC{\begin{corollary}}
\def\EC{\end{corollary}}
\def\BE{\begin{example}}
\def\EE{\end{example}}
\def\BD{\begin{definition}}
\def\ED{\end{definition}}
\def\BR{\begin{remark}}
\def\ER{\end{remark}}
\def\BAS{\begin{assumption}}
\def\EAS{\end{assumption}}
\def\BI{\begin{itemize}}
\def\EI{\end{itemize}}

\def\BMP{\begin{minipage}{9.5cm}}
\def\EMP{\end{minipage}}
\def\MPT{\begin{minipage}{11.5cm}}
\def\EPT{\end{minipage}}

\def\H{\mathbb{H}}
\def\R{\mathbb{R}}

\newcommand*\diff{\mathop{}\!\mathrm{d}}

\newcommand*\samethanks[1][\value{footnote}]{\footnotemark[#1]}
\newcommand{\NNorm}[2]{\left\Vert {#1} \right\Vert_{#2}}
\newcommand{\AV}[1]{\left\vert {#1} \right\vert}

\title{Optimality conditions for mathematical programs \\ with orthogonality type constraints}
\author{
S. L\"ammel
\thanks{
Department of Mathematics, Chemnitz University of Technology,
Reichenhainer Str. 41, 09126
Chemnitz, Germany; e-mail: sebastian.laemmel@mathematik.tu-chemnitz.de, vladimir.shikhman@mathematik.tu-chemnitz.de
(corresponding author).
 } \and V. Shikhman\samethanks[1]
}

\begin{document}
\thispagestyle{empty}
\maketitle
\vspace{-5ex}
\abstract{ We consider the class of mathematical programs with orthogonality type constraints (MPOC). Orthogonality type constraints appear by reformulating the sparsity constraint via auxiliary binary variables and relaxing them afterwards. For MPOC a necessary optimality condition in terms of T-stationarity is stated. The justification of T-stationarity is threefold. First, it allows to capture the global structure of MPOC in terms of Morse theory, i.\,e.~deformation and cell-attachment results are established. For that, nondegeneracy for the T-stationary points is introduced and shown to hold at a generic  MPOC. Second, we prove that Karush-Kuhn-Tucker points of the Scholtes-type regularization converge to T-stationary points of MPOC. This is done under the MPOC-tailored linear independence constraint qualification (LICQ), which turns out to be a generic property too. Third, we show that T-stationarity applied to the relaxation of sparsity constrained nonlinear optimization (SCNO) naturally leads to its M-stationary points. Moreover, we argue that all T-stationary points of this relaxation become degenerate.
}

\vspace{2ex}
{\bf Keywords: orthogonality type constraints, T-stationarity, nondegeneracy, genericity, Morse theory, Scholtes-type regularization, sparsity constraint}

\vspace{2ex}
{\bf MSC-classification: 90C26}

\section{Introduction}
\label{sec:intro}
We consider the following mathematical program with orthogonality type constraints:
\[
\mbox{MPOC}: \quad
\min_{x} \,\, f(x)\quad \mbox{s.\,t.} \quad x \in M
\]
with
\[
    M=\left\{x \in\R^n \left\vert\;\begin{array}{l}
    h_i(x)=0,i\in I, 
    g_j(x)\ge 0, j \in J,\\ 
    F_{1,m}(x) \cdot F_{2,m}(x)=0, F_{2,m}(x)\ge 0, m=1,\ldots,k 
    \end{array} \right.\right\},
\]
  where $f \in C^2(\R^n,\R)$, $ h \in C^2(\R^n, \R^{|I|})$, $ g \in C^2(\R^n, \R^{|J|})$,
  $F_1, F_2 \in C^2(\R^n,\R^k)$.
Orthogonality type constraints appeared initially in \cite{burdakov:2016}. There, the so-called sparsity constrained nonlinear optimization problem (SCNO) has been reformulated as a mixed-integer program by means of auxiliary variables. By relaxing the binary constraints, the authors obtained a particular instance of MPOC, see also Section 5 for details. 

Our aim here is rather to study the class of MPOC in its generality. Crucial for our approach is the introduction of T-stationarity for MPOC. T-stationarity is not only a necessary optimality condition, see Theorem \ref{prop:necess}, but, more importantly, it is relevant from the topological point of view. In order to become specific on this matter, let us briefly overview our findings on T-stationarity and its justification: 
\begin{itemize}
    \item[(1)] The global structure of MPOC is captured by T-stationary points. Namely, two basic theorems within the scope of Morse theory are shown, cf. \cite{jongen:2000}. The deformation result states that outside the set of T-stationary points the topology of MPOC lower level sets remains unchanged if the level varies, see Theorem \ref{thm:def}. The cell-attachment result describes the topological changes of MPOC level sets if passing a nondegenerate T-stationary point, see Theorem \ref{thm:cell-a}. More precisely, the dimension of the cell to be attached corresponds to the T-index. The latter has not only the usual quadratic part, but  includes also the linear part associated with biactive orthogonality type constraints. We point out that this is quite usual in structured nonsmooth optimization, cf. \cite{jongen:2009} for mathematical programs with complementarity constraints (MPCC) and \cite{dorsch:2012} for mathematical programs with vanishing constraints (MPVC).
    \item[(2)] The Scholtes-type regularization of orthogonality type constraints naturally leads to T-stationary points. 
    Namely, Karush-Kuhn-Tucker points of the Scholtes-type regularization converge to T-stationary points of MPOC, see Theorem \ref{thm:regul}. This is shown to hold under the MPOC-tailored linear independence constraint qualification (LICQ). This result is in accordance with the observation made also for other classes of nonsmooth optimization problems. Often, Karush-Kuhn-Tucker points of the Scholtes-type regularization approach topologically relevant stationary points, cf. \cite{mehlitz:2019} and \cite{shikhman:2021} for mathematical programs with switching constraints (MPSC).
    \item[(3)] Finally, we apply the notion of T-stationarity to the relaxation of SCNO from \cite{burdakov:2016}. It turns out that the T-stationary points of this relaxation correspond to M-stationary points of SCNO. Note that M-stationarity is known to be topologically relevant for SCNO in the sense of Morse theory, see \cite{laemmel:2019}. Although this appealing relation between both stationarity concepts,
    the relaxation of SCNO introduces intrinsic degeneracies. Namely, all its T-stationary points are shown to become degenerate, even if the corresponding M-stationary points of SCNO were originally nondegenerate, see Theorem \ref{thm:deg}. This is the price to pay for rewriting the sparsity constraint by means of auxiliary variables. 
\end{itemize}

The article is organized as follows. In Section 2 we introduce the notion of T-stationarity and provide auxiliary results which will be used later.
Section 3 contains the exposition of Morse theory for MPOC including the proofs of the deformation and cell-attachment theorems.
Section 4 is devoted to the Scholtes-type regularization of MPOC and its relation to T-stationarity. In Section 5 we apply T-stationarity to the relaxation of SCNO and conclude on its degeneracy.

Our notation is standard. The $n$-dimensional Euclidean space is denoted by $\R^n$ with the coordinate vectors $e_i$, $i=1, \ldots,n$. The $n$-dimensional vector with all components equal to one is denoted by $e$. We denote the set of nonnegative numbers by $\mathbb{H}$. The solution set of the basic orthogonality type relation is 
$
  \mathbb{T} = \left\{\left. \left(a, b\right) \in \R^2 \, \right| \, a \cdot b =0, b \geq 0\right\}$.
Given a differentiable function $F:\R^n \longrightarrow \R^m$, $DF$ denotes its Jacobian matrix.
Given a differentiable function $f:\R^n \longrightarrow \R$, $Df$ denotes the row vector of partial derivatives.

\section{Optimality conditions}

Given $\bar x \in M$, we define the following index sets:
\[
J_0(\bar x)=\left\{j \in J \; \vert\; g_j(\bar x)=0\right\},
\]
    \[a_{00}\left(\bar x\right)=\left\{m\in\left\{1,\ldots,k\right\}\left\vert F_{1,m}\left(\bar x\right)=0,F_{2,m}\left(\bar x\right)=0\right.\right\},\]
     \[a_{01}\left(\bar x\right)=\left\{m\in\left\{1,\ldots,k\right\}\left\vert F_{1,m}\left(\bar x\right)=0,F_{2,m}\left(\bar x\right)>0\right.\right\},\] \[a_{10}\left(\bar x\right)=\left\{m\in\left\{1,\ldots,k\right\}\left\vert F_{1,m}\left(\bar x\right)\ne 0,F_{2,m}\left(\bar x\right)=0\right.\right\}.\]
The index set $J_0(\bar x)$ corresponds to the active inequality constraints and $\alpha_{00}({\bar x})$ to the bi-active orthogonality type constraints at $\bar x$.
Without loss of generality, we assume throughout the whole article that
at the particular point of interest $\bar x \in M$ it holds:
\[
J_0(\bar x)=\{1, \ldots, |J_0(\bar x)|\},
\]
\[
\alpha_{01}(\bar x)=\{1, \ldots, |\alpha_{01}(\bar x)|\}, \quad \alpha_{10}(\bar x)=\{|\alpha_{01}(\bar x)|+1, \ldots, |\alpha_{01}(\bar x)|+|\alpha_{10}(\bar x)|\},
\]
\[
a_{00}\left(\bar x\right)=\left\{ 
|a_{01}(\bar x)|+|a_{10}(\bar x)|+1,\ldots,k\right\}.
\]
Furthermore, we put:
\[
s=|I|+|a_{01}(\bar x)|+|a_{10}(\bar x)|,\quad q=s+|J_0(\bar x)|, \quad p=n-q-2|a_{00}\left(\bar x\right)|.
\]

Let us start by stating the MPOC-tailored linear independence constraint qualification, which turns out to be the crucial assumption for all results to follow.

\begin{definition}[LICQ]
We say that a feasible point $\bar x \in M$ of MPOC satisfies the linear independence constraint qualification (LICQ)  if the following vectors are linearly independent:
\[
\begin{array}{l}
D h_i(\bar x), i \in I,
D g_j(\bar x), j \in J_0(\bar x), \\
D F_{1,m}(\bar x), m \in a_{01}\left(\bar x\right)\cup a_{00}\left(\bar x\right),
D F_{2,m}(\bar x), m \in a_{10}\left(\bar x\right)\cup a_{00}\left(\bar x\right).
\end{array}
\]
\end{definition}
The assumption of LICQ is justified in the sense, that it generically holds on the MPOC feasible set. In order to formulate this assertion in mathematically precise terms, the space $C^2(\R^n,\R)$ will be topologized by means of the strong (or Whitney-) $C^2$-topology, denoted by $C^2_s$, cf. \cite{hirsch:1976}, \cite{jongen:2000}. The $C^2_s$-topology is generated by allowing
perturbations of the functions and their derivatives up to second order which are controlled by means of continuous
positive functions. The product space $C^2(\R^n, \R^l) \cong C^2(\R^n, \R) \times \cdots \times C^2(\R^n, \R)$ will be topologized with the corresponding product topology.

\begin{proposition}[Genericity of LICQ]
\label{prop:licq-gen}
Let $\mathcal{F}\subset C^2\left(\R^n,\R^{n+|I|+|J|+k}\right)$
denote the subset of MPOC defining functions for which LICQ holds at all feasible points. Then, $\mathcal{F}$ is
$C_s^2$-open and -dense.
\end{proposition}
\proof
We consider the feasible set of a disjunctive optimization problem:
\[
M_{DISJ}=\left\{x\in \R^n \left\vert \;\begin{array}{l}
    h_i(x)=0,i\in I, 
    g_j(x)\ge 0, j \in J,\\ 
    \max\left\{F_{1,m}(x), F_{2,m}(x)\right\}\ge 0, m=1,\ldots,k 
    \end{array}  \right.\right\}.
\]
The set $M$ is clearly a subset of $M_{DISJ}\left\lbrack h,g,F_1,F_2\right\rbrack$.
In \cite{jongen:1997}, it was shown that LICQ holds generically on $M_{DISJ}$. Since the notions of LICQ are the same for disjunctive optimization problems and MPOC, the assertion follows immediately. 
\qed

Next, we introduce a stationarity concept for MPOC which appears to be relevant for our topologically motivated studies. 

\begin{definition}[T-stationary point]
\label{def:t-stat}
A feasible point $\bar x \in M$ is called T-stationary for MPOC if there exists multipliers $$\bar \lambda_i, i\in I,
\bar \mu_j, j\in J_0(\bar x),\bar \sigma_{1,m_{01}}, m_{01} \in a_{01}\left(\bar x\right),\bar \sigma_{2,m_{10}}, m_{10} \in a_{10}\left(\bar x\right),\bar \varrho_{1,m_{00}},\bar \varrho_{2,m_{00}}, m_{00} \in a_{00}\left(\bar x\right),$$ such that the following conditions hold:
\begin{equation}
   \label{eq:tstat-1} 
   \begin{array}{rcl}
   D f(\bar x)&=& \sum\limits_{i\in I}\bar \lambda_i D h_i(\bar x)+
    \sum\limits_{j \in J_0(\bar x)}\bar \mu_j D g_j(\bar x) \\ \\
    && +\sum\limits_{m_{01} \in a_{01}\left(\bar x\right)} \bar \sigma_{1,m_{01}}DF_{1,m_{01}}\left(\bar x\right)
    +\sum\limits_{m_{10} \in a_{10}\left(\bar x\right)} \bar \sigma_{2,m_{10}}DF_{2,m_{10}}\left(\bar x\right) \\ \\
    &&+\sum\limits_{m_{00} \in a_{00}\left(\bar x\right)} \left(\bar \varrho_{1,m_{00}}DF_{1,m_{00}}\left(\bar x\right)+\bar \varrho_{2,m_{00}}DF_{2,m_{00}}\left(\bar x\right)\right), \end{array}
\end{equation}
\begin{equation}
   \label{eq:tstat-2} \bar \mu_j \ge 0 \mbox{ for all } j\in J_0\left(\bar x\right),
\end{equation}
\begin{equation}
   \label{eq:tstat-3} \bar \varrho_{1,m_{00}}=0 \mbox{ or }\bar \varrho_{2,m_{00}}\le 0 \mbox{ for all } m_{00} \in a_{00}\left(\bar x\right).
\end{equation}

\end{definition}
It turns out that a necessary optimality condition for MPOC can be stated in terms of T-stationarity. To show this, we first locally describe the MPOC feasible set under LICQ. This is done by an appropriate change of coordinates.

\begin{lemma}[Local structure]
\label{lem:local}
Suppose that LICQ holds at $\bar x \in M$. Then $M$ admits a local $C^2$-coordinate system of $\R^n$ at $\bar x$, i.\,e.~there exists a $C^2$-diffeomorphism $\Phi: U \rightarrow V$ with open 
$\R^n$-neighborhoods $U$ and $V$ of $\bar x$ and $0$, respectively, such that it holds:
\begin{itemize}
    \item[(i)] $\Phi(\bar x)=0$,
    \item[(ii)] $\Phi\left(M \cap U\right) = \left(\{0\}^{s} \times \H^{\left|J_0(\bar x)\right|} \times \mathbb{T}^{|a_{00}\left(\bar x\right)|} \times \R^{p}\right) \cap V$.
\end{itemize} 
\end{lemma}
\proof
We choose vectors $\xi_1,\ldots,\xi_{p}\in \R^n$,
which together with the vectors 
\[
\begin{array}{l}
D h_i(\bar x), i \in I, 
D g_j(\bar x), j \in J_0(\bar x), \\
D F_{1,m}(\bar x), m \in a_{01}\left(\bar x\right)\cup a_{00}\left(\bar x\right), 
D F_{2,m}(\bar x), m \in a_{10}\left(\bar x\right)\cup a_{00}\left(\bar x\right)
\end{array}
\]
form a basis for $\R^n$.
We put:
\[
\begin{array}{lcl}
    y_i &=&h_i\left( x\right) \mbox{ for } i \in I,  \\
    y_{|I|+m_{01}}&=& F_{1,m_{01}}\left( x\right) \mbox{ for } m_{01} \in a_{01}(\bar x),\\
     y_{|I|+m_{10}}&=& F_{2,m_{10}}\left( x\right) \mbox{ for } m_{10} \in a_{10}(\bar x),\\
    y_{s+j}&=&g_j\left( x\right) \mbox{ for } j \in J_0(\bar x), \\
    y_{q+2(m_{00}-|a_{01}(\bar x)|-|a_{10}(\bar x)|)-1}&=& F_{1,m_{00}}\left( x\right) \mbox{ for } m_{00} \in a_{00}(\bar x),\\
    y_{q+2(m_{00}-|a_{01}(\bar x)|-|a_{10}(\bar x)|)}&=& F_{2,m_{00}}\left( x\right) \mbox{ for } m_{00} \in a_{00}(\bar x),\\
 y_{n-p+r}&=&\xi_{r}^T\left(x-\bar x\right) \mbox{ for } r=1,\ldots,p.
 \end{array}
\]
We write for short 
\begin{equation}
  \label{eq:stddiff}
y=\Phi(x).  
\end{equation}
By definition it holds $\Phi(x) \in C^2\left(\R^n,\R^n\right)$ and $\Phi(\bar x)=0$. Due to LICQ, the Jacobian matrix $D\Phi(\bar x)$ is nonsingular.
Hence, by means of the inverse function theorem, there exist open neighborhoods
$U$ of $\bar x$ and $V$ of $0$ such that $\Phi:U\rightarrow V$ is a $C^2$-diffeomorphism. Moreover, we can guarantee that $J_0(x) \subset J_0(\bar x)$ and $a_{00}(x)\subset a_{00}(\bar x)$ by shrinking $U$ if necessary. Thus, property (ii) follows directly from the definition of $\Phi$.\qed

\begin{theorem}[Necessary optimality condition]
\label{prop:necess}
Let $\bar x$ be a local minimizer of MPOC satisfying LICQ, then $\bar x$ is a T-stationary point.
\end{theorem}

\proof Due to LICQ at $\bar x$, we may write $Df(\bar x)$ as a linear combination of the 
basis vectors given in the proof of Lemma \ref{lem:local}:
\begin{equation}
   \label{eq:tstat-help} 
   \begin{array}{rcl}
   D f(\bar x)&=& \sum\limits_{i\in I}\bar \lambda_i D h_i(\bar x)+
    \sum\limits_{j \in J_0(\bar x)}\bar \mu_j D g_j(\bar x) \\ \\
    && +\sum\limits_{m_{01} \in a_{01}\left(\bar x\right)} \bar \sigma_{1,m_{01}}DF_{1,m_{01}}\left(\bar x\right)
    +\sum\limits_{m_{10} \in a_{10}\left(\bar x\right)} \bar \sigma_{2,m_{10}}DF_{2,m_{10}}\left(\bar x\right) \\ \\
    &&+\sum\limits_{m_{00} \in a_{00}\left(\bar x\right)} \left(\bar \varrho_{1,m_{00}}DF_{1,m_{00}}\left(\bar x\right)+\bar \varrho_{2,m_{00}}DF_{2,m_{00}}\left(\bar x\right)\right) +\sum\limits_{r=1}^{p} \bar \nu_r \xi_r(\bar x). 
    \end{array}
\end{equation}
It is straightforward to see that the unique multipliers in (\ref{eq:tstat-help}) are derivatives of the objective function $f \circ \Phi^{-1}$ in new coordinates given by the diffeomorphism $\Phi$ from (\ref{eq:stddiff}), cf. \cite[Lemma 2.2.1]{jongen:2004}. Due to Lemma \ref{lem:local}, $\bar x$ is a local minimizer of $f$ on the MPOC feasible set if and only if $0$ is a local minimizer of $f \circ \Phi^{-1}$ on 
$\{0\}^{s} \times \H^{\left|J_0(\bar x)\right|} \times \mathbb{T}^{|a_{00}\left(\bar x\right)|} \times \R^{p}$. The latter provides in particular:
\[
   \frac{\partial(f \circ \Phi^{-1})(0)}{\partial y_{s+j}} \geq 0, j \in J_0(\bar x), \quad 
   \frac{\partial(f \circ \Phi^{-1})(0)}{\partial y_{q+2(m_{00}-|a_{01}(\bar x)|-|a_{10}(\bar x)|)-1}}=0, m_{00} \in \alpha_{00}(\bar x), 
   \]
   \[
   \frac{\partial(f \circ \Phi^{-1})(0)}{\partial y_{n-p+r}}=0, r=1, \ldots,p.
\]
or, equivalently,
\[
   \bar \mu_j \geq 0, j \in J_0(\bar x), \quad 
   \bar \varrho_{1,m_{00}} =0, m_{00} \in \alpha_{00}(\bar x), \quad \bar \nu_r=0, r=1, \ldots,p.
\]
Thus, (\ref{eq:tstat-1})-(\ref{eq:tstat-3}) from the definition of T-stationarity follow. \qed

As we have seen in the proof of Theorem \ref{prop:necess}, the multipliers $(\bar \lambda, \bar \mu, \bar \sigma, \bar \varrho)$ of a T-stationary point $\bar x$ are the corresponding partial derivatives of the objective function $f\circ\Phi^{-1}$ in new coordinates given by the diffeomorphism $\Phi$ from (\ref{eq:stddiff}).
Let us now also interpret  the second-order derivatives of the objective function $f \circ \Phi^{-1}$ in new coordinates. For that, given a T-stationary point $\bar x \in M$ with the multipliers $\left(\bar \lambda, \bar \mu, \bar \sigma, \bar \varrho\right)$ we define the Lagrange function:
\[
        \begin{array}{rcl}
        L(x)&=&\displaystyle f(x)-
         \sum_{i \in I} \bar \lambda_i h_i( x)-\sum_{j \in J_0(\bar x)} \bar \mu_j g_j( x) \\ \\
        &&-\sum\limits_{m_{01} \in a_{01}\left(\bar x\right)}  \bar \sigma_{1,m_{01}}DF_{1,m_{01}}\left(\bar x\right)-\sum\limits_{m_{10} \in a_{10}\left(\bar x\right)}  \bar \sigma_{2,m_{10}}DF_{2,m_{10}}\left(\bar x\right) \\ \\
    &&-\sum\limits_{m_{00} \in a_{00}\left(\bar x\right)} \left( \bar \varrho_{1,m_{00}}DF_{1,m_{00}}\left(\bar x\right)- \bar \varrho_{2,m_{00}}DF_{2,m_{00}}\left(\bar x\right)\right). \end{array}
\]
 We further consider the corresponding local part of the feasible set around $\bar x$:
\[
  \begin{array}{ll}
      M(\bar x)= \{ x \in \R^n \,\,| & h_i(x) =0, i \in I,
    g_j(x) =0, j \in J_0(\bar x),  \\ &
       F_{1,m_{01}}(x) = 0, m_{01} \in \alpha_{01}(\bar x), F_{2,m_{10}}(x) = 0, m_{10} \in \alpha_{01}(\bar x) , \\ &
    F_{1,m_{00}}(x) = 0, F_{2,m_{00}}(x) = 0, m_{00} \in \alpha_{00}(\bar x)
       \}.
     \end{array}
\]
Obviously, $M(\bar x) \subset M$ in a sufficiently small neighborhood of $\bar x$. Moreover, in the case where LICQ holds at $\bar x$, $M(\bar x)$ is locally a
$p$-dimensional $C^2$-manifold. The tangent space of $M(\bar x)$ at $\bar x$ is
\[
    \begin{array}{ll}
      T_{\bar x}M(\bar x) = \{ \xi \in \R^n \,\,| &
        D h_i(\bar x) \, \xi=0, i \in I, D g_j(\bar x) \, \xi =0, j \in J_0(\bar x), \\ &
        D F_{1,m_{01}}(\bar x) \, \xi = 0, m_{01} \in \alpha_{01}(\bar x), 
        D F_{2,m_{10}}(\bar x) \, \xi= 0, m_{10} \in \alpha_{10}(\bar x), 
     \\ &
    D F_{1,m_{00}}(\bar x) \, \xi = 0, D F_{2,m_{00}}(\bar x) \, \xi = 0, m_{00} \in \alpha_{00}(\bar x)
         \}.
     \end{array}
\]
 It is easy to see that the Hessian of $f \circ \Phi^{-1}$ with respect to the last $p$ coordinates corresponds to the restriction of the Lagrange function's Hessian $D^2 L(\bar x)$ on the respective tangent space $T_{\bar x}M(\bar x)$, cf. \cite[Lemma 2.2.10]{jongen:2004}. This observation motivates to introduce the notion of nondegeneracy for T-stationary points. 

\begin{definition}[Nondegenerate T-stationary point]
\label{def:nondeg}
A T-stationary point $\bar x \in M$ of MPOC with multipliers $(\bar \lambda, \bar \mu, \bar \sigma, \bar \varrho)$ is called nondegenerate if the following conditions are satisfied:
\begin{itemize}
    \item []ND1: LICQ holds at $\bar x$,
    \item []ND2: the strict complementarity (SC) holds for active inequality constraints, i. e. $\bar \mu_j>0$ for all $j \in J_0\left(\bar x\right)$,
    \item []ND3: the multipliers corresponding to biactive orthogonality type constraints do not vanish, i.\,e. $\bar\varrho_{1,m_{00}}\ne 0$ and $\bar\varrho_{2,m_{00}}< 0$ for all $m_{00}\in a_{00}\left(\bar x\right)$,
    \item []ND4: the matrix $D^2 L(\bar x)\restriction_{T_{\bar x} M(\bar x)}$ is nonsingular.
\end{itemize}  
Otherwise, we call $\bar x$ degenerate.
\end{definition}
The assumption of nondegeneracy is justified in the sense, that it generically holds at all T-stationary points of MPOC.

\begin{proposition}
[Nondegeneracy is generic]
\label{prop:generic}
Let $\mathcal{F}\subset C^2\left(\R^n,\R^{n+|I|+|J|+k}\right)$ be the subset of MPOC defining functions for which each 
T-stationary point is nondegenerate. Then, $\mathcal{F}$ is $C^2_s$-open and -dense.
\end{proposition}

\proof    Let us fix an index set $ J_0 \subset J$ of active inequality constraints, an index subset $T_0 \subset J_0$ of those active inequality constraints, three pairwise disjoint sets $a_{01},a_{10}, a_{00}\subset \left\{1,\ldots,k\right\}$ with $ a_{01}\cup a_{10}\cup a_{00}=\left\{1,\ldots,k\right\}$ representing the three active index sets. Furthermore, we fix some subsets $\alpha_1, \alpha_2 \subset a_{00}$ of the biactive index set, and a number $r \in \mathbb{N}$ standing for the rank.
For this choice we consider the set  $M_{J_0,T_0,a_{01},a_{10},a_{00},\alpha_1,\alpha_2,r}$ of $x \in \R^n$ such that the following conditions are satisfied:
\begin{itemize}
    \item[] (m1) $h_i(x)=0$ for all $i \in I$ and $g_{j}(x)=0$ for all $j \in J_0$,
    \item[] (m2) $F_{1,m_1}(x)=0$ for all $m_1 \in \left(a_{00} \cup a_{01}\right)$ and $F_{2,m_2}(x)=0$ for all $m_2 \in \left(a_{00} \cup a_{10}\right)$,
    \item[] (m3) 
    \[
    \begin{array}{rcl}
   D f(x)&=& \sum\limits_{i\in I} \lambda_i D h_i( x)+
    \sum\limits_{J_0\backslash T_0} \mu_j D g_j(x) \\ \\
    && +\sum\limits_{m_{01} \in a_{01}}  \sigma_{1,m_{01}}DF_{1,m_{01}}\left( x\right)
    +\sum\limits_{m_{10} \in a_{10}}  \sigma_{2,m_{10}}DF_{2,m_{10}}\left( x\right) \\ \\
    &&+\sum\limits_{m_{00} \in a_{00}\backslash \alpha_1}  \varrho_{1,m_{00}}DF_{1,m_{00}}\left( x\right)+ \sum\limits_{m_{00} \in a_{00}\backslash \alpha_2}  \varrho_{2,m_{00}}DF_{2,m_{00}}\left( x\right), \end{array}
    \]
    \item[] (m4) the matrix $D^2 L(x)\restriction_{T_{x} M(x)}$ has rank $r$.
\end{itemize}
Note that (m1) refers to equality and active inequality constraints, (m2) to constraints regarding the functions $F_1$ and $F_2$, while (m3) describes violation of ND2 and ND3. Furthermore, (m4) describes the violation of ND4.
Now, it suffices to show $M_{J_0,T_0,a_{01},a_{10},a_{00},\alpha_1,\alpha_2,r}$ is generically empty whenever one of the sets $T_0,\alpha_1$ or $\alpha_2$ is nonempty or the rank $r$ in (m4) is not full, i.\,e. $r < \mbox{dim}\left(T_{x} M(x)\right)$.
In fact, the available degrees of freedom
of the variables involved in each $M_{J_0,T_0,a_{01},a_{10},a_{00},\alpha_1,\alpha_2,r}$ are $n$. The loss of freedom caused by (m1) is $\left|I\right|+\left|J_0\right|$, and
the loss of freedom caused by (m2) is $\left|a_{01}\right|+\left|a_{10}\right|+2\left|a_{00}\right|$.
Due to Proposition \ref{prop:licq-gen}, LICQ holds generically at any feasible $x$, i.\,e.~(ND1) is fulfilled. Suppose that the sets $T_0,\alpha_1$, and $\alpha_2$ are empty, then (m3) causes a loss of freedom of $n-\left|I\right|-\left|
J_0\right|-(\left|a_{01}\right|+\left|a_{10}\right|+2\left|a_{00}\right|)$. Hence, the total loss of freedom is $n$.
We conclude that a
further degeneracy, i.\,e. $T_0 \not = \emptyset$, $\alpha_1 \not = \emptyset$, $\alpha_2 \not = \emptyset$  or $r < \mbox{dim}\left(T_{x} M(x)\right)$, would imply that the total available degrees of freedom $n$ are exceeded. By virtue of the jet
transversality theorem from \cite{jongen:2000}, generically the sets  $M_{J_0,T_0,a_{01},a_{10},a_{00},\alpha_1,\alpha_2,r}$ must be empty.
For the openness result, we argue in a standard way. Locally, T-stationarity can be written
via stable equations. Then, the implicit function theorem for Banach spaces can be applied to
follow T-stationary points with respect to (local) $C^2$-perturbations of defining functions. Finally,
a standard globalization procedure exploiting the specific properties of the strong $C^2_s$-topology can be used to construct a (global) $C^2_s$-neighborhood of problem data for which the nondegeneracy
property is stable, cf. \cite{jongen:2000}. \qed

With a nondegenerate T-stationary point it is convenient to associate its T-index. 

\begin{definition}[T-index]
Let $\bar x \in M$ be a nondegenerate T-stationary point of MPOC with unique multipliers $\left(\bar \lambda,\bar \mu, \bar \sigma,\bar \varrho\right)$. The  number of negative eigenvalues of the matrix $D^2 L(\bar x)\restriction_{T_{\bar x} M(\bar x)}$ is called its quadratic index ($QI$). The number of negative $\bar \varrho_{2,m}$ with $m \in a_{00}\left(\bar x\right)$ equals $\left|a_{00}\left(\bar x\right)\right|$ and is called the biactive index ($BI$) of $\bar x$. We define the T-index ($TI$) as the sum of both, i.\,e. $TI=QI+BI$.
\end{definition}
The T-index encodes the local structure of MPOC around the corresponding T-stationary point.

\begin{proposition}[Morse Lemma for MPOC]
\label{prop:morse}
Suppose that $\bar x$ is a nondegenerate T-stationary point of MPOC with quadratic index $QI$ and biactive index $BI$. Then, there exist neighborhoods $U$ and $V$ of $\bar x$ and $0$, respectively, and a local $C^1$-coordinate system $\Psi: U \rightarrow V$ of $\R^n$ around $\bar x$ such that:

\begin{equation}
\label{eq:normal}
    f\circ \Psi^{-1}(y)= f(\bar x) +
    \sum\limits_{j \in J_0(\bar x)} y_{s+j} +
    \sum\limits_{m=1}^{\left|a_{00}\left(\bar x\right)\right|}\left( \pm y_{q+2m-1}-y_{q+2m}\right) + \sum\limits_{r=1}^p\pm y_{n-p+r}^2,  
\end{equation}
where $y \in \{0\}^{s} \times \H^{\left|J_0(\bar x)\right|} \times \mathbb{T}^{|a_{00}\left(\bar x\right)|} \times \R^{p}$. Moreover, there are exactly $QI$ negative squares in (\ref{eq:normal}).
\end{proposition}


\proof
Without loss of generality, we may assume $f(\bar x)=0$.
By using $\Phi$ from (\ref{eq:stddiff}), we put $\bar f := f \circ \Phi^{-1}$ on the set 
$\left(\{0\}^{s} \times \H^{\left|J_0(\bar x)\right|} \times \mathbb{T}^{|a_{00}\left(\bar x\right)|} \times \R^{p}\right)\cap V$.
At the origin we have with respect to the new $y$-coordinates, see Lemma \ref{lem:local}:
\begin{itemize}
    \item [(i)] $\displaystyle \frac{\partial \bar f}{\partial y_{s+j}} > 0$ for $j \in J_0(\bar x)$,
    \item [(ii)] $\displaystyle \frac{\partial \bar f}{\partial y_{q+2m-1}}\ne 0,\displaystyle \frac{\partial \bar f}{\partial y_{q+2m}}<0$ for $m\in\left\{1,\ldots,\left|a_{00}\left(\bar x\right)\right|\right\}$,
    \item [(iii)] $\displaystyle \frac{\partial \bar f}{\partial y_{n-p+r}}=0$ for $r\in \left\{1,\ldots,p\right\}$ and the matrix $\displaystyle \left(\frac{\partial^2 \bar f}{\partial y_{n-p+r_1} \partial y_{n-p+r_2}}\right)_{r_1,r_2 \in \left\{1,\ldots,p\right\}}$ is non\-singular.
\end{itemize}
We denote $\bar f$ by $f$ again. Under the following coordinate transformations the set $$\{0\}^{s} \times \H^{\left|J_0(\bar x)\right|} \times \mathbb{T}^{|a_{00}\left(\bar x\right)|} \times \R^{p}$$
will be equivariantly transformed in itself. We put $y=\left(Y_{n-r},Y^r\right)$, where 
\[
Y_{n-p}=\left(y_1, \ldots, y_{n-p}\right), \quad
Y^{p}=\left(y_{n-p+1},\ldots,y_n\right).
\]

It holds:
\[
\begin{array}{lcl}
f\left(Y_{n-p},Y^{p}\right)&=&
\displaystyle \int_0^1 \frac{d}{dt}f\left(tY_{n-p},Y^{p}\right)\diff t+f\left(0,Y^p\right)\\ \\
&=&\displaystyle \sum\limits_{r=1}^{n-p-s} y_{s+r} \cdot d_{s+r}(y) + f\left(0,Y^p\right),
\end{array}
\]
where
\[
   d_{s+r}(y)=\displaystyle \int_0^1 \frac{\partial f}{\partial y_{s+r}}\left(tY_{n-p},Y^{p}\right)\diff t,\quad r=1,\ldots, n-p-s.
\]
Note that 
 $d_{s+r} \in C^1$ for $r=1,\ldots, n-p-s.$
Due to (iii), we may apply the standard Morse Lemma on the $C^2$-function $f\left(0,Y^{p}\right)$ without affecting the first $Y_{n-p}$ coordinates, see e.\,g. \cite{jongen:2000}. The corresponding coordinate transformation is of class $C^1$. Denoting the transformed functions again by $f$ and $d_{s+r}$, we obtain
\[
f(y)= \displaystyle \sum\limits_{r=1}^{n-p-s} y_{s+r} \cdot d_{s+r}(y)+ \sum\limits_{r \in R}\pm y_{n-p+r}^2.
\]

Furthermore, (i) and (ii) provide that
\[
\begin{array}{rcl}
   d_{s+j}(0)&=&\displaystyle \frac{\partial f}{\partial y_{s+j}}\left(0\right) > 0,\quad j \in J_0(\bar x),
\\ \\
 d_{q+2m-1}(0)&=&\displaystyle \frac{\partial f}{\partial y_{q+2m-1}}\left(0\right)\ne 0,\quad m\in\left\{1,\ldots,\left|a_{00}\left(\bar x\right)\right|\right\},
 \\ \\
  d_{q+2m}(0)&=&\displaystyle \frac{\partial f}{\partial y_{q+2m}}\left(0\right)< 0,\quad m\in\left\{1,\ldots,\left|a_{00}\left(\bar x\right)\right|\right\}.
\end{array}
\]
Hence, we may take
\[
\begin{array}{l}
y_{s+j} \cdot \left|d_{s+j}(y)\right|, \quad j \in J_0(\bar x), \\
y_{q+2m-1} \cdot d_{q+2m-1}(y), \quad  m\in\left\{1,\ldots,\left|a_{00}(\left(\bar x\right)\right|\right\}, \\
y_{q+2m} \cdot \left|d_{q+2m}(y)\right|, \quad  m\in\left\{1,\ldots,\left|a_{00}(\left(\bar x\right)\right|\right\}, \\
y_{n-p+r}, \quad r \in R
\end{array}
\]
as new local $C^1$-coordinates by a straightforward application of the inverse function theorem. Denoting the transformed function again by $f$, we obtain (\ref{eq:normal}). Here, the coordinate transformation $\Psi$
is understood as the composition of all previous ones.
\qed
As a by-product we elaborate on how the T-index can be used to characterize nondegenerate local minimizers of MPOC. This result can be seen as a sufficient optimality condition for generic MPOCs.

\begin{corollary}[Sufficient optimality condition]
\label{lem:min-index}
Let $\bar x \in M$ be a nondegenerate T-stationary point. Then, $\bar x$ is a local
minimizer of MPOC if and only if its T-index vanishes.
\end{corollary}

\proof Let $\bar x$ be a nondegenerate T-stationary point for MPOC.
The application of Morse Lemma from Proposition \ref{prop:morse} says that there exist
neighborhoods $U$ and $V$ of $\bar x$ and $0$, respectively, and a local $C^1$-coordinate system $\Psi: U \rightarrow V$ of $\R^n$ around $\bar x$ such that (\ref{eq:normal}) holds. Therefore, $\bar x$ is a local minimizer for MPOC if and only if $0$ is a local minimizer of $ f\circ \Psi^{-1}$ on the set $\{0_s\} \times \mathbb{H}^{|J_0(\bar x)|} \times \mathbb{T}^{|\alpha_{00}(\bar x)|} \times \R^{p}$.
If the T-index vanishes, we have $\alpha_{00}(\bar x)=\emptyset$ and $QI=0$, and (\ref{eq:normal}) reads as
 \begin{equation}
    \label{eq:normalatMI0}
    f \circ \Psi^{-1} \left( 0_s, y_{s+1}, \ldots, y_n \right) =
    f(\bar x)+
    \sum_{j=1}^{|J_0(\bar x)|} y_{s+j} 
    +\sum_{r=1}^{p} y^2_{n-p+r},
 \end{equation}
where $y \in \{0_s\} \times \mathbb{H}^{|J_0(\bar x)|} \times \R^{p}$. Thus, $0$ is a local minimizer for (\ref{eq:normalatMI0}).
Vice versa, if $0$ is a local minimizer for (\ref{eq:normalatMI0}), then obviously $\alpha_{00}(\bar x)=\emptyset$ and $QI=0$, hence, the T-index of $\bar x$ vanishes.
\qed 

Corollary \ref{lem:min-index} says that at a nondegenerate local  minimizer $\bar x \in M$ of MPOC with multipliers $(\bar \lambda, \bar \mu, \bar \sigma,\bar \varrho)$ we have LICQ and SC, but also the following conditions are satisfied:
\begin{itemize}
    \item[(a)]  the bi-active orthogonality type constraints are absent, i.\,e. $\alpha_{00}(\bar x)=\emptyset$,
    \item[(b)]  the second-order sufficient condition (SOSC) is fulfilled, i.\,e.~the matrix $D^2 L(\bar x) \restriction_{T_{\bar x}M(\bar x)}$ is positive definite.
\end{itemize}
 Additionally, the application of Proposition \ref{prop:generic} provides that all minimizers of MPOC are generically nondegenerate. In particular, their structure is relatively simple with respect to the orthogonality constraints due to (a). The importance of nondegeneracy for local minimizers becomes clear if we turn our attention to stability issues. Example \ref{ex:stab} provides that the violation of condition (a) at a local minimizer may lead to instability w.r.t. arbitrarily small $C^2$-perturbations of the MPOC defining functions.
 
 \begin{example}[Instability]
     \label{ex:stab}  
     Consider the MPOC:
\begin{equation}
  \label{eq:exinst1}
  \min \,\, x_1^2 + x_2^2 \quad \mbox{s.t.} \quad x_1 \cdot x_2 = 0, x_2 \geq 0. 
\end{equation}
Obviously, $(0,0)$ uniquely solves (\ref{eq:exinst1}). It is a T-stationary point for (\ref{eq:exinst1}) with both bi-active multipliers vanishing, hence, it is degenerate.
Consider the following perturbation of (\ref{eq:exinst1}) with respect to parameter $\varepsilon$:
\begin{equation}
  \label{eq:ex1_pert}
  \min \,\,  \left(x_1+\varepsilon\right)^2 + \left(x_2-\varepsilon\right)^2 \quad \mbox{s.t.} \quad x_1 \cdot x_2 = 0, x_2 \geq 0.
\end{equation}
It is easy to see that $(0,0)$, $(0,\varepsilon)$ and $(-\varepsilon,0)$ are T-stationary points for (\ref{eq:ex1_pert}), where the parameter $\varepsilon>0$ is sufficiently small, but arbitrary. More than that, both $(0,\varepsilon)$ and $(-\varepsilon,0)$ solve (\ref{eq:ex1_pert}). 
This observation suggests that $(0,0)$ is unstable as the solution for (\ref{eq:exinst1}). \qed
 \end{example}




\section{Morse theory}

We study the topological properties of MPOC lower level sets
\[
M_a=\left\{x\in M \left\vert f(x)\le a\right.\right\},
\]
where $a\in \R$ is varying. For that, we define intermediate sets for $a < b$:
\[
M^b_a=\left\{x\in M \left\vert a\le f(x) \le b\right. \right\}.
\]

Outside the set of T-stationary points, the topology of the MPOC lower level sets remains unchanged. This is referred to as deformation result within the scope of Morse theory, cf. \cite{jongen:2000}. 

\begin{theorem}[Deformation for MPOC]
\label{thm:def}
Let $M^b_a$ be compact and LICQ be fulfilled at all points $x \in M_a^b$. Then, if $M_a^b$ contains no T-stationary points for MPOC, then $M_a$ is homeomorphic to $M_b$.
\end{theorem}

\proof
For all $x\in M_a^b$ there exist due to LICQ multipliers 
$$\lambda_i(x), i\in I,
\mu_j(x), j\in J_0(x),\sigma_{1,m_{01}}(x), m_{01} \in a_{01}\left(x\right),\sigma_{2,m_{10}}(x), m_{10} \in a_{10}\left(x\right),\varrho_{1,m_{00}}(x),$$ 
$$\varrho_{2,m_{00}}(x), m_{00} \in a_{00}\left( x\right),\nu_r(x), r=1,\ldots,p,$$ such that
\[
\begin{array}{rl}
   D f( x)=& \sum\limits_{i\in I} \lambda_i(x) D h_i( x)+
    \sum\limits_{j \in J_0(x)}\mu_j(x) D g_j(x) \\ \\
    & +\sum\limits_{m_{01} \in a_{01}\left( x\right)}  \sigma_{1,m_{01}}(x) DF_{1,m_{01}}\left(x\right)
    +\sum\limits_{m_{10} \in a_{10}\left( x\right)} \sigma_{2,m_{10}}(x) DF_{2,m_{10}}\left( x\right) \\ \\
    &+\sum\limits_{m_{00} \in a_{00}\left( x\right)} \left( \varrho_{1,m_{00}}(x) DF_{1,m_{00}}\left( x\right)+ \varrho_{2,m_{00}}(x) DF_{2,m_{00}}\left( x\right)\right) + \sum\limits_{r=1}^p \nu_r(x)\xi_r
    , \end{array}
\]

where the vectors $\xi_r, r\in R$ are chosen as in Lemma \ref{lem:local}.
Next, we set:
\[
\begin{array}{rcl}
A&=& \left\{x \in M_a^b \left\vert \mbox{there exists } r\in \{1, \ldots,p\} \mbox{ such that } \nu_r(x)\ne 0\right. \right\},\\
B&=& \left\{x \in M_a^b \left\vert \mbox{there exists } j\in J_0(x) \mbox{ such that } \mu_j(x)< 0\right. \right\},\\
C&=&\left\{x \in M_a^b \left\vert \mbox{there exists } m_{00} \in a_{00}(x) \mbox{ such that } \varrho_{1,m_{00}}(x)\ne 0 \mbox{ and } \varrho_{2,m_{00}}(x)>0 \right.\right\}.
\end{array}
\]
For $\bar x \in M_a^b$ we get $\bar x \in A \cup B \cup C$, since it is not T-Stationary for MPOC.
The proof consists of a local argument an its globalization. First, we show the local argument, i.\,e., for each
$\bar x \in M_a^b$ there exist a neighborhood $U_{\bar x}$ of $\bar x, t_{\bar x}>0$, and a mapping
\[
\Psi_{\bar x}:\left\{
\begin{array}{rcl}
\left[0,t_{\bar x}\right) \times \left(M^b\cap U_{\bar x}\right)&\longrightarrow &M,\\
\left(t,x\right)&\mapsto&\Psi_{\bar x}\left(t,x\right),
\end{array}\right.
\]
such that
\begin{itemize}
    \item [(i)] $\Psi_{\bar x}\left(t,M^b\cap U_{\bar x}\right)\subset M^{b-t}$ for all $t \in \left[0,t_{\bar x}\right)$,
    \item [(ii)] $\Psi_{\bar x}\left(t_1+t_2,\cdot\right)=\Psi_{\bar x}\left(t_1,\Psi_{\bar x}\left(t_2,\cdot\right)\right)$ for all $t_1,t_2 \in \left[0,t_{\bar x}\right)$ with $t_1+t_2 \in \left[0,t_{\bar x}\right)$,
    \item [(iii)] if $\bar x \in A\cup B$, then $\Psi_{\bar x}\left(\cdot,\cdot\right)$ is a $C^1$-flow corresponding to a $C^1$-vector field $F_{\bar x}$,
    \item[(iv)] if $\bar x \in C$, then $\Psi_{\bar x}\left(\cdot,\cdot\right)$ is a Lipschitz flow.
\end{itemize}

Note that the level sets of $f$ are mapped locally onto the level sets of $f\circ \Phi^{-1}$, where $\Phi$ is the the diffeomorphism from (\ref{eq:stddiff}), see Lemma \ref{lem:local}. We consider  $f\circ \Phi^{-1}$ and denote it by $f$ again.
Thus, we have $\bar x=0$ and $f$ is given on the feasible set
$\{0\}^{s} \times \H^{\left|J_0(\bar x)\right|} \times \mathbb{T}^{|a_{00}\left(\bar x\right)|} \times \R^{p}$.

Case (a): $\bar x \in A$\\
It follows that there exists $r\in \{1,\ldots,p\}$ with $\displaystyle \frac{\partial f}{\partial x_r}\left(\bar x\right)\ne 0$. We define a local $C^1$-vector field $F_{\bar x}$ as
\[
F_{\bar x}(x_1,\ldots,x_r,\ldots,x_n)=\left(0,\ldots,\displaystyle-\frac{\partial f}{\partial x_r}\left(x\right)\cdot\left(\displaystyle\frac{\partial f}{\partial x_r}\left(x\right)\right)^{-2},\ldots,0\right)^T,
\]
which -- after respective inverse changes of local coordinates -- induces the flow $\Psi_{\bar x}$ fitting the local argument.

Case (b): $\bar x \in B$\\
It follows that there exists $j\in J_0(\bar x)$ with $\displaystyle \frac{\partial f}{\partial x_{s+j}}\left(\bar x\right)<0$. By means of a local $C^1$-coordinate transformation in the $(s+j)$-th coordinate on $\H$, leaving the other coordinates unchanged, we obtain locally for $f$
\[
f\left(x_1,\ldots,x_{s+j},\ldots,x_n\right)=-x_{s+j}+f\left(x_1,\ldots,\bar x_{s+j},\ldots,x_n\right)
\]
and define
\[
F_{\bar x}(x_1,\ldots,x_{s+j},\ldots,x_n)=\left(0,\ldots,1,\ldots,0\right)^T,
\]
which -- after respective inverse changes of local coordinates -- induces the flow $\Psi_{\bar x}$ fitting the local argument.

Case (c): $\bar x \in C$\\
It follows that there exists $m\in \left\{1,\ldots,\left|a_{00}(\bar x)\right|\right\}$ with $\displaystyle \frac{\partial f}{\partial x_{q+2m-1}}\left(\bar x\right)\ne 0$ and
$\displaystyle \frac{\partial f}{\partial x_{q+2m}}\left(\bar x\right) > 0$.
We define $\delta=\mbox{sign}\left(\displaystyle \frac{\partial f}{\partial x_{q+2m-1}}\left(\bar x\right)\right)$.
Analogously to the proof of Proposition \ref{prop:morse}, we obtain the following representation for $f$ in $C^1$-coordinates:
\[
f(x)=\delta \cdot x_{q+2m-1}+x_{q+2m}+f\left(x_1,\ldots,\bar x_{q+2m-1},\bar x_{q+2m},\ldots, x_n\right).
\]
We define the mapping $\Psi_{\bar x}$ locally as
\[
\begin{array}{l}
\Psi_{\bar x}\left(t,x_1,\ldots,x_{q+2m-1}, x_{q+2m},\ldots, x_n\right)\\=
\left(x_1,\ldots, x_{q+2m-1}-\delta \cdot \max \left\{0,t-x_{q+2m}\right\} ,\max\{0,x_{q+2m}-t\},\ldots, x_n\right)^T.
\end{array}
\]
Obviously, the latter also fits the local argument.

Next, we globalize the local argument.
Therefore, consider the open covering $\left\{U_x\left\vert\, x\in C \right.\right\}\cup\left\{U_{\bar x} \left\vert\, \bar x \in M_a^b\backslash\left\{U_x\left\vert\, x \in C\right.\right\} \right.\right\}$. 
The sets $U_{\bar x},\,\bar x\in M_a^b\backslash\left\{U_x\left\vert\, x \in C\right.\right\}$ can be chosen smaller and if necessary disjoint with $C$ by using continuity arguments.
Since $M_a^b$ is compact, we get a finite subcovering $\left\{U_{x^{\ell}}\left\vert\,x^{\ell}\in C\right.\right\}\cup \left\{U_{\bar x_k} \left\vert \bar x_k \in M_a^b\backslash\left\{U_x\left\vert\, x \in C\right.\right\} \right.\right\}$ of $M_a^b$. Let $\left\{\phi_k\right\}$ be a $C^{\infty}$-partition of unity subordinate to this subcovering. We define the $C^1$-vector field
$F=\sum\limits_k \phi_kF_{x_k}$, which induces a flow on $\left\{U_{\bar x_k} \left\vert \bar x_k \in M_a^b \right.\right\}$ using $F_{x_k}$ from cases (a) and (b). The latter induces a flow $\Psi$ on
$\left\{U_{\bar x_k} \left\vert \bar x_k \in M_a^b\backslash\left\{U_x\left\vert\, x \in C\right.\right\} \right.\right\}$. In each nonempty overlapping region $U_{x^{\ell}}\cap U_{\bar x_k}, x^{\ell}\in C, \bar x_k \in 
 M_a^b\backslash\left\{U_x\left\vert\, x \in C\right.\right\}$ the flow $\Psi_{x^{\ell}}$ induces exactly the vector field $F$. Hence, local trajectories can be glued together on $M_a^b$. We denote the resulting mapping by $\Psi$ again. Furthermore, the level of $f$ reduces at least by
 \[
 \displaystyle\frac
 {\min\left\{t_{x^{\ell}},t_{\bar x_k}\left\vert\, x^{\ell} \in C, \bar x_k \in M_a^b\backslash\left\{U_x\left\vert\, x \in C\right.\right\}\right.\right\}}{2}>0.
 \]
 Thus, for $x\in M_a^b$ we obtain a unique $t_a(x)>0$ with $\Psi\left(t_a(x),x\right)\in M^a$.
 The mapping $t_a:x\rightarrow t_a(x)$ is indeed Lipschitz. Finally,
we define $r:[0,1]\times M^b\rightarrow M^b$ as
\[
r\left(\tau,x\right)=
\left\{
\begin{array}{ll}
x&\mbox{for }x \in M^a,\tau \in [0,1]\\
\Psi\left(\tau\cdot t_a(x),x\right)&\mbox{for }x \in M_a^b,\tau \in [0,1].
\end{array}
\right.
\]
This mapping provides that $M_a$ is a strong deformation retract of $M^b$.
\qed

Let us now turn our attention to the topological changes of lower level sets when passing a T-stationary level. Traditionally, they are described by means of the so-called cell-attachment.
We first consider a special case of cell-attachment. For that, let $N^\epsilon$ denote the lower level set of a special linear function on $\H^\ell \times \mathbb{T}^m$, i.\,e.
\[
 N^{\epsilon} = \left\{ \left(u,v\right) \in \H^\ell \times \mathbb{T}^m \,\left|\, \sum\limits_{i=1}^{\ell} u_i + \sum\limits_{j=1}^{m} \left(\pm v_{2j-1}-v_{2j}\right)\leq \epsilon \right. \right\},
\]
where $\epsilon \in \R$, and the integers $\ell$ and $m$ are nonnegative. 

\begin{lemma}[Normal Morse data]
\label{lem:cat}
For any $\epsilon > 0$ the set $N^\epsilon$ is homotopy-equivalent to $N^{-\epsilon}$ with a $m$-dimensional cell attached.
\end{lemma}
\proof
 The lower level set $N^{-\epsilon}$ of a special linear function on $\H^\ell \times \mathbb{T}^m$ is given by
\[
 N^{-\epsilon} = \left\{ \left(u,v\right) \in \H^\ell \times \mathbb{T}^m \,\left|\, \sum\limits_{i=1}^{\ell} u_i + \sum\limits_{j=1}^{m} \left(\pm v_{2j-1}-v_{2j}\right)\leq -\epsilon \right. \right\},
\]
with $\epsilon >0$. This is homotopy-equivalent to the set
\[
\bar  N^{-\epsilon} = \left\{ v \in \mathbb{T}^m \,\left|\, \sum\limits_{j=1}^{m} \left(\pm v_{2j-1}-v_{2j}\right)\leq -\epsilon \right. \right\},
\]
using the homotopy
\[
\left((u,v),t\right)) \mapsto \left((1-t) \cdot u,v\right), \quad t \in [0,1].
\]
 The lower level set $\bar N^{-\epsilon}$ itself is homotopy-equivalent to the set
 \[
 \hat N^{-\epsilon} =\left\{ v \in \R^{2m} \left\vert\;
 \begin{array}{l}
v_{2j-1}\cdot v_{2j}=0 \mbox { and } v_{2j-1},v_{2j}\ge  0 \mbox{ for all } j=1,\ldots, m,\\ 
\sum\limits_{j=1}^{m} \left(-v_{2j-1}-v_{2j}\right)\leq -\epsilon  \end{array}\right. \right\}
 \]
by using the homotopy
\[
\left(v,t\right) \mapsto \left( (1-t)\cdot v_{2j-1}\mp t \cdot \left|v_{2j-1}\right|, v_{2j}, j=1, \ldots, m\right), \quad t \in [0,1].
\]
Note that $\hat N^{-\epsilon}$ is a lower level set of the following mathematical program with complemantary constaints (MPCC):
\begin{equation}
\label{eq:MPCC}
\min_{v \in R^{2m}} \,\, \sum\limits_{j=1}^{m} \left(-v_{2j-1}-v_{2j}\right)\quad \mbox{s.\,t.} \quad v \in \hat M,
\end{equation}
where
\[
   \hat M=\left\{v \in \R^{2m} \left\vert\;
 v_{2j-1}\cdot v_{2j}=0 \mbox { and } v_{2j-1}, v_{2j}\ge 0 \mbox{ for all } j=1,\ldots, m\right.\right\}.
\]
The cell-attachment for mathematical programs with complemantary constraints was examined in 
\cite{jongen:2009}. Note that the origin is the only C-stationary point of (\ref{eq:MPCC}), moreover, it is nondegenerate with C-index equal to $m$. According to the results in \cite{jongen:2009},
$\hat N^{\epsilon}$ is then homotopy-equivalent to $\hat N^{-\epsilon}$ with an $m$-dimensional cell attached. 
Since both $\hat N^{\epsilon}$ and $N^{\epsilon}$ are contractible given the homotopy
\[
\left(v,t\right)) \mapsto \left((1-t)\cdot v\right), \quad t \in [0,1],
\]
the assertion follows.
\qed

Now, we are ready to describe in general how the topology of the MPOC lower level sets changes if passing a T-stationary level.

\begin{theorem}[Cell-Attachment for MPOC]
\label{thm:cell-a}
Let $M_a^b$ be compact and suppose that it contains exactly one nondegenerate T-stationary point $\bar x$ with T-index equal to $t$. 
If $a<f\left(\bar x \right) <b$,
then $M^b$ is homotopy-equivalent to $M^a$ with a $t$-cell attached.
\end{theorem}
\proof
Theorem \ref{thm:def} allows deformations up to an arbitrarily small neighborhood of the M-stationary point $\bar x$. In such a neighborhood, we may assume without loss of generality that $\bar x=0$ and $f$ has the following form as from Proposition \ref{prop:morse}:
\begin{equation}
\label{eq:c-att}
    f\circ \Psi^{-1}(y)= f(\bar x) +
    \sum\limits_{j \in J_0(\bar x)} y_{s+j} +
    \sum\limits_{m=1}^{\left|a_{00}\left(\bar x\right)\right|}\left( \pm y_{q+2m-1}-y_{q+2m}\right) + \sum\limits_{r=1}^p\pm y_{n-p+r}^2,  
\end{equation}
where $y \in \{0\}^{s} \times \H^{\left|J_0(\bar x)\right|} \times \mathbb{T}^{|a_{00}\left(\bar x\right)|} \times \R^{p}$ and the number of negative squares in (\ref{eq:c-att}) equals $QI$.

In terms of \cite{goresky:1988} the set 
$\{0\}^{s} \times \H^{\left|J_0(\bar x)\right|} \times \mathbb{T}^{|a_{00}\left(\bar x\right)|} \times \R^{p}$ can be interpreted as the product of the tangential part $\{0\}^{s} \times \R^{p}$ and the normal part $\H^{\left|J_0(\bar x)\right|} \times \mathbb{T}^{|a_{00}\left(\bar x\right)|}$. The cell-attachment along the tangential part is standard.  Analogously to the case of nonlinear programming, one $QI$-dimensional cell has to be attached on $\{0\}^{s} \times \R^{p}$.  The cell-attachment along the normal part can be deduced from Lemma \ref{lem:cat}: a $BI$-dimensional cell has to be attached.
Finally, we apply Theorem 3.7 from Part I in \cite{goresky:1988}, which says that the local Morse data is the product of tangential and normal Morse data. Hence, the dimensions of the attached cells add together. Here, we have then to attach a $(QI+BI)$-dimensional cell on $\{0\}^{s} \times \H^{\left|J_0(\bar x)\right|} \times \mathbb{T}^{|a_{00}\left(\bar x\right)|} \times \R^{p}$. The dimension corresponds to the T-index $t$. \qed

A global interpretation of our results is typical for the Morse theory, cf. \cite{jongen:2000}. Suppose that the MPOC feasible
set is compact and connected, that 
LICQ holds at all feasible points, and that all T-stationary points are nondegenerate with pairwise
different functional values. Then, passing a level corresponding to a local minimizer,
a connected component of the lower level set is created. Different components can
only be connected by attaching 1-cells. This shows the existence of at least $(k-1)$ T-stationary
points with T-index equal to one, where $k$ is the number of local minimizers. In the context of global optimization the latter result usually is referred to as mountain pass. For nonlinear programming it is well known that Karush-Kuhn-Tucker points with quadratic index equal to one naturally appear along with local minimizers, see e.\,g. \cite{floudas:2005}. For MPOC, however, not only T-stationary points with quadratic index equal to one, but also with biactive index equal to one may become relevant, see next Example \ref{ex:global}.

\begin{example}[Saddle point]
  \label{ex:global}
  We consider the following MPOC:
\begin{equation}
  \label{eq:st1}
   \min_{x_1,x_2}\,\, \left(x_1+1\right)^2 + \left(x_2-1\right)^2 \quad \mbox{s.\,t.} \quad 
   x_1 \cdot x_2=0, x_2 \geq 0.
\end{equation}
Obviously, $(-1,0)$ and $(0,1)$ are nondegenerate minimizers for (\ref{eq:exinst1}). Hence, there should exist an additional T-stationary point with T-index one. 
This nondegenerate T-stationary point is $(0,0)$. In fact, the orthogonality type constraint is biactive at $(0,0)$. Moreover, the multiplier $\varrho_1=2$ corresponding to $x_1$ is positive and the multiplier $\varrho_2=-2$ corresponding to $x_2$ is negative. Thus, its quadratic index vanishes and its biactive index equals one, i.\,e. $QI=0$, $BI=1$.\qed
\end{example}

\section{Scholtes-type regularization}

Let us consider the Scholtes-type regularization of MPOC with the parameter $t >0$:
\[
\mbox{MPOC}_t: \quad
\min_{x} \,\, f(x)\quad \mbox{s.\,t.} \quad x \in M_{t},
\]
where
\begin{equation}
    M_{t}=\left\{x \in\R^n \left\vert\;\begin{array}{l}
    h_i(x)=0,i\in I, 
    g_j(x)\ge 0, j \in J,\\ 
    -t \le F_{1,m}(x) \cdot F_{2,m}(x)\le t,
    F_{2,m}(x)\ge 0, m=1,\ldots,k 
    \end{array} \right.\right\}.
\end{equation}
We now prove that the Karush-Kuhn-Tucker points of the regularized problem $\mbox{MPOC}_t$ converge to a T-stationary point of MPOC if $t$ goes to zero. 
Similar results have been shown in \cite{branda:2018} for cardinality constrained nonlinear optimization problems (CCOP) and in \cite{kanzow:2019} for mathematical programs with switching constraints (MPSC). The proof will follow along the lines of those papers.

\begin{theorem}[Convergence for Scholtes-type regularization]
\label{thm:regul}
Let $\left\{t^{\ell}\right\}_{{\ell}\in \N}$  be a sequence of positive regularization parameters converging to zero. Suppose for each ${\ell}\in \N$ there exists a Karush-Kuhn-Tucker point $x^{\ell}\in M_{t^\ell}$ of 
$\mbox{MPOC}_{t^\ell}$, and let the resulting sequence $\left\{x^{\ell}\right\}_{\ell\in \N}$ converge to $\bar x$. If LICQ holds at $\bar x \in M$, then it is a T-stationary point of MPOC.
\end{theorem}
\proof
Since $x^{\ell}$ is a Karush-Kuhn-Tucker point for $\mbox{MPOC}_{t^{\ell}}$, there exist multipliers
\[
\lambda_i^{\ell}, i \in I, \mu_j^{\ell}, j \in J_0\left(x^{\ell}\right),
\eta_m^{\ell},\eta_m^{{\ell},\geq},\eta_m^{{\ell},\leq}, m=1, \ldots,k,
\]
such that it holds:
\begin{itemize}
    \item []S1:\[
    \begin{array}{rcl}
   D f\left(x^{\ell}\right)&=& \sum\limits_{i\in I}\lambda^{\ell}_i D h_i\left(x^{\ell}\right)+
    \sum\limits_{j \in J}\mu^{\ell}_j D g_j\left(x^{\ell}\right)+
    \sum\limits_{m=1}^k\eta^{\ell}_m D F_{2,m}\left(x^{\ell}\right)\\ \\
    && +\sum\limits_{m=1}^k \eta^{{\ell},\geq}\left(DF_{1,m}\left(x^{\ell}\right)\cdot F_{2,m}\left(x^{\ell}\right)+DF_{2,m}\left(x^{\ell}\right)\cdot F_{1,m}\left(x^{\ell}\right)\right)\\ \\
    && -\sum\limits_{m=1}^k \eta^{{\ell},\leq}\left(DF_{1,m}\left(x^{\ell}\right)\cdot F_{2,m}\left(x^{\ell}\right)+DF_{2,m}\left(x^{\ell}\right)\cdot F_{1,m}\left(x^{\ell}\right)\right), \end{array}\]
    \item []S2: $\mu^\ell_j \cdot g_j\left(x^\ell\right)=0 \mbox{ and } \mu^{\ell}_j \ge 0 \mbox{ for all } j\in J$,
    \item []S3: $ \eta^{\ell}_m \cdot F_{2,m}\left(x^{\ell}\right)=0 \mbox{ and }  \eta^{\ell}_m\ge 0 \mbox{ for all } m=1,\ldots,k $,
    \item []S4a: $\eta^{\ell,\geq}_m\cdot \left(F_{1,m}\left(x^{\ell}\right) \cdot F_{2,m}\left(x^{\ell}\right)+t^{\ell}\right)=0 \mbox{ and } \eta^{\ell,\geq}_m\ge 0 \mbox{ for all } m=1,\ldots,k $,
    \item []S4b: $\eta^{\ell,\leq}_m\cdot \left(t^{\ell}-F_{1,m}\left(x^{\ell}\right) \cdot F_{2,m}\left(x^{\ell}\right)\right)=0 \mbox{ and } \eta^{\ell,\leq}_m\ge 0 \mbox{ for all } m=1,\ldots,k $.
\end{itemize}
Next, we define for $\ell \in \N$ new multipliers:
\begin{itemize}
    \item []$\sigma_{1,m}^{\ell}=\left\{\begin{array}{ll}
   \eta_m^{\ell,\ge}\cdot F_{2,m}\left(x^{\ell}\right)-\eta_m^{\ell,\le}\cdot F_{2,m}\left(x^{\ell}\right),& \mbox{for }m \in a_{01}\left(\bar x\right),\\
    0,& \mbox{else},\end{array}\right.$
     \item []$\sigma_{2,m}^{\ell}=\left\{\begin{array}{ll}
    \eta_m^{\ell}+\eta_m^{\ell,\ge}\cdot F_{1,m}\left(x^{\ell}\right)-\eta_m^{\ell,\le}\cdot F_{1,m}\left(x^{\ell}\right),& \mbox{for }m \in a_{10}\left(\bar x\right),\\
    0,& \mbox{else},\end{array}\right.$
    \item []$\varrho_{1,m}^{\ell}=\left\{\begin{array}{ll}
   \eta_m^{\ell,\ge}\cdot F_{2,m}\left(x^{\ell}\right)-\eta_m^{\ell,\le}\cdot F_{2,m}\left(x^{\ell}\right),& \mbox{for }m \in a_{00}\left(\bar x\right),\\
    0,& \mbox{else},\end{array}\right.$
    \item []$\varrho_{2,m}^{\ell}=\left\{\begin{array}{ll}
      \eta_m^{\ell}+\eta_m^{\ell,\ge}\cdot F_{1,m}\left(x^{\ell}\right)-\eta_m^{\ell,\le}\cdot F_{1,m}\left(x^{\ell}\right),& \mbox{for }m \in a_{00}\left(\bar x\right),\\
    0,& \mbox{else}.\end{array}\right.$
\end{itemize}
Note that for sufficiently large $\ell \in \N$ we have $\eta^{\ell}_{m_{01}}=0$, $m_{01} \in a_{01}({ \bar x})$. This follows by applying continuity arguments to the definition of $a_{01}({ \bar x})$ and S3.
Thus, we can replace S1 by the following equation:
 \begin{equation}
     \label{eq:reg1}
   \begin{array}{rcl}
   D f( x^{\ell})&=&\sum\limits_{i\in I}\lambda^{\ell}_i D h_i\left(x^{\ell}\right)+
    \sum\limits_{j \in J}\mu^{\ell}_j D g_j\left(x^{\ell}\right)\\ \\
   &&+\sum\limits_{m_{01}\in a_{01}({ \bar x})} \sigma^{\ell}_{1,m_{01}} D F_{1,m_{01}}({ x^{\ell}})+\sum\limits_{m_{10}\in a_{10}({ \bar x})}\sigma^{\ell}_{2,m_{10}} D F_{2,m_{10}}({ x^{\ell}})
      \\ \\
    &&+\sum\limits_{m_{00}\in a_{00}({\bar x})}
   \varrho^{\ell}_{1,m_{00}} DF_{1,m_{00}}\left(x^{\ell}\right)
   +\varrho^{\ell}_{2,m_{00}} DF_{2,m_{00}}\left(x^{\ell}\right)\\ \\
    && +\sum\limits_{m_{01}\in a_{01}({ \bar x})} \left(\eta^{{\ell},\ge}-\eta^{{\ell},\le}\right)\cdot F_{1,m_{01}}\left(x^{\ell}\right)\cdot D F_{2,m_{01}}\left(x^{\ell}\right)\\ \\&&+\sum\limits_{m_{10}\in a_{10}({ \bar x})}\left(\eta^{{\ell},\ge}-\eta^{{\ell},\le}\right)\cdot F_{2,m_{10}}\left(x^{\ell}\right)\cdot D F_{1,m_{10}}\left(x^{\ell}\right).
    \end{array}
 \end{equation}

Our first claim is that the sequence of multipliers $\left(\lambda^{\ell},\mu^{\ell},\eta^{{\ell},\ge}-\eta^{{\ell},\le},\sigma^\ell,\varrho^\ell\right)$ is bounded on the set $a_{01}\left(\bar x\right)$.
If not, we define for $\ell \in \N$:
\[
\left(\tilde\lambda^{\ell},\tilde\mu^{\ell},\tilde\eta^{{\ell},\ne},\tilde\sigma^{\ell},\tilde\varrho^{\ell}\right)=
\frac{\left(\lambda^{\ell},\mu^{\ell},\eta^{{\ell},\ge}-\eta^{{\ell},\le},\sigma^{\ell},\varrho^{\ell}\right)}
{\left\|\left(\lambda^{\ell},\mu^{\ell},\eta^{{\ell},\ge}-\eta^{{\ell},\le},\sigma^{\ell},\varrho^{\ell}\right)\right\|_2}.
\]
The resulting sequence  is bounded and, thus, without loss of generality converges to a vector $\left(\tilde\lambda,\tilde\mu,\tilde\eta^{\ne},\tilde\sigma,\tilde\varrho\right)$ (after choosing a suitable subsequence if necessary). 
We divide (\ref{eq:reg1}) by the norm of the multipliers and take the limit $\ell \rightarrow \infty$:
\begin{equation}
    \label{eq:reg-t1}
    \begin{array}{rcl}
   0&=&\sum\limits_{i\in I}\tilde\lambda_i D h_i(\bar x)+
    \sum\limits_{j \in J}\tilde\mu_j D g_j(\bar x)\\ \\
   &&+\sum\limits_{m_{01}\in a_{01}({ \bar x})} \tilde\sigma_{1,m_{01}} D F_{1,m_{01}}({ \bar x})+\sum\limits_{m_{10}\in a_{10}({ \bar x})}\tilde\sigma_{2,m_{10}} D F_{2,m_{10}}({ \bar x})\\ \\
    &&+\sum\limits_{m_{00}\in a_{00}({ \bar x})}
   \tilde\varrho_{1,m_{00}} DF_{1,m_{00}}(\bar x)
   +\tilde\varrho_{2,m_{00}} DF_{2,m_{00}}(\bar x).
    \end{array}
    \end{equation}
Due to LICQ at $\bar x$, the multipliers $\tilde \lambda, \tilde \mu,\tilde \sigma,\tilde \varrho$ vanish.
Hence, there exists an index $m_{01}\in a_{01}\left(\bar x\right)$ such that $\tilde \eta_{m_{01}}^{\ne}\ne 0$. Moreover, it holds:
\[
\sigma^{\ell}_{1,m_{01}} = \left(\eta_{m_{01}}^{{\ell},\ge}-\eta_{m_{01}}^{{\ell},\le}\right)\cdot F_{2,m_{01}}\left(x^{\ell}\right).
\]
Dividing by $\left\|\left(\lambda^{\ell},\mu^{\ell},\eta^{{\ell},\ge}-\eta^{{\ell},\le},\sigma^{\ell},\varrho^{\ell}\right)\right\|_2$ and taking the limit leads to
\[
 \tilde\sigma_{1,m_{01}}=\tilde\eta_{m_{01}}^{\ne}F_{2,m_{01}}\left(\bar x\right)\ne 0.
\]
This contradicts the fact that $\tilde\sigma$ vanishes, and we are done. 

Our second claim is that the sequence
$\left\{\left(\lambda^{\ell},\mu^{\ell},\eta^{{\ell},\ge}-\eta^{{\ell},\le},\sigma^{\ell},\varrho^{\ell}\right)\right\}_{\ell \in \N}$ is bounded on the set $\left\{m_{10}\right\}$ with an arbitrary $m_{10} \in a_{10}\left(\bar x\right)$ or for $\ell \in \N$ large enough it holds $F_{2,m_{10}}\left(x^{\ell}\right)=0$.
In order to show this, we suppose the sequence to be not bounded on $\left\{m_{10}\right\}$ and conclude -- similarly to the case of $a_{01}\left(\bar x\right)$ -- that the multipliers $\tilde \lambda, \tilde \mu,\tilde \sigma,\tilde \varrho$ vanish, but $\tilde\eta^{\ne}_{m_{10}}\ne 0$ holds. 
If, additionally, $F_{2,m_{10}}\left(x^{\ell_0}\right)\ne 0$ along some subsequence $\ell_0$, we have by using S3:
\[
\sigma^{\ell_0}_{2,m_{10}}=\underbrace{\eta^{\ell_0}_{m_{10}}}_{=0}+\left(\eta_{m_{10}}^{{\ell_0},\ge}-\eta_{m_{10}}^{{\ell_0},\le}\right)\cdot F_{1,m_{10}}\left(x^{\ell_0}\right).
\]
Dividing by $\left\|\left(\lambda^{\ell},\mu^{\ell},\eta^{{\ell},\ge}-\eta^{{\ell},\le},\sigma^{\ell},\varrho^{\ell}\right)\right\|_2$ and taking the limit leads to
\[
\tilde\sigma_{2,m_{10}}= \tilde\eta_{m_{10}}^{\not=}F_{1,m_{10}}\left(\bar x\right)\ne 0.
\]
This contradicts the fact that $\tilde\sigma$ vanishes, and we are done again. 

Now, we are ready to take the limit $\ell \rightarrow \infty$ in (\ref{eq:reg1}). Our first claim provides that for every $m_{01}\in a_{01}\left(\bar x\right)$ we have:
\[
\left(\eta^{{\ell},\ge}-\eta^{{\ell},\le}\right)\cdot F_{1,m_{01}}\left(x^{\ell}\right)\cdot D F_{2,m_{01}}\left(x^{\ell}\right) \rightarrow 0 \quad \mbox{if } \ell \rightarrow \infty.
\]
Our second claim provides that
for every $m_{10}\in a_{10}\left(\bar x\right)$ we have:
\[
\left(\eta^{{\ell},\ge}-\eta^{{\ell},\le}\right)\cdot F_{2,m_{10}}\left(x^{\ell}\right)\cdot D F_{1,m_{10}}\left(x^{\ell}\right) \rightarrow 0 \quad \mbox{if } \ell \rightarrow \infty.
\]
Due to LICQ at $\bar x$, the sequence $\left\{\left(\lambda^{\ell},\mu^{\ell},\sigma^{\ell},\varrho^{\ell}\right)\right\}_{\ell \in \N}$   converges without loss of generality to a vector
$\left(\bar \lambda,\bar \mu,\bar \sigma,\bar \varrho\right)$ (after choosing a suitable subsequence if necessary).
Thus, (\ref{eq:reg1}) provides in the limit:
\[
     \begin{array}{rcl}
      D f\left(\bar x\right)&=&\sum\limits_{i\in I}\bar\lambda_i D h_i(\bar x)+
    \sum\limits_{j \in J}\bar\mu_j D g_j(\bar x)\\ \\
   &&+\sum\limits_{m_{01}\in a_{01}({ \bar x})} \bar\sigma_{1,m_{01}} D F_{1,m_{01}}({ \bar x})+\sum\limits_{m_{10}\in a_{10}({ \bar x})}\bar\sigma_{2,m_{10}} D F_{2,m_{10}}({ \bar x})\\ \\
    &&+\sum\limits_{m_{00}\in a_{00}({ \bar x})}
   \bar\varrho_{1,m_{00}} DF_{1,m}(\bar x)
   +\bar\varrho_{2,m_{00}} DF_{2,m_{00}}(\bar x).
    \end{array}\]
Due to continuity reasons, it holds $J_0\left(x^{\ell}\right)\subset J_0\left(\bar x\right)$, hence,
\[
   \bar \mu_j \cdot g_j\left(\bar x\right)=0 \mbox{ and } \bar \mu_j \ge 0 \mbox{ for all } j\in J.
\]
It remains to prove that $ \bar \varrho_{1,m_{00}}=0 \mbox{ or }\bar \varrho_{2,m_{00}}\le 0 \mbox{ for all } m_{00} \in a_{00}\left(\bar x\right)$. If $\bar \varrho_{1,m_{00}}\ne 0$, then $\varrho^\ell_{1,m_{00}}\ne 0$ for $\ell \in \N$ large enough, and, due to their definition it holds:
\[
 \eta_{m_{00}}^{\ell,\ge}-\eta_{m_{00}}^{\ell,\le} \ne 0 \quad \mbox{and} \quad  F_{2,{m_{00}}}\left(x^{\ell}\right)> 0.
\]
Hence, $\eta^{\ell}_{m_{00}}=0$ for those $\ell \in \N$ by S3.
Moreover, by adding S4a and S4b we obtain:
\begin{equation}
\label{eq:reg2}
\left(\eta^{\ell,\geq}_{m_{00}}-\eta^{\ell,\leq}_{m_{00}}\right)\cdot F_{1,{m_{00}}}\left(x^{\ell}\right) \cdot F_{2,{m_{00}}}\left(x^{\ell}\right)+t^{\ell}\cdot \left(\eta^{\ell,\geq}_{m_{00}}+\eta^{\ell,\leq}_{m_{00}}\right)=0.
\end{equation}
By simplifying the latter, we get:
\begin{equation}
\label{eq:reg3}
  \left(\eta^{\ell,\geq}_{m_{00}}-\eta^{\ell,\leq}_{m_{00}}\right)\cdot F_{1,{m_{00}}}\left(x^{\ell}\right)=\frac{-t^{\ell}\cdot \left(\eta^{\ell,\geq}_{m_{00}}+\eta^{\ell,\leq}_{m_{00}}\right)}{F_{2,{m_{00}}}{\left(x^{\ell}\right)}}\le 0.
\end{equation}
Thus, we have:
\[
\varrho_{2,{m_{00}}}^{\ell}= \eta^{\ell}_{m_{00}}+ \left(\eta^{\ell,\geq}_{m_{00}}-\eta^{\ell,\leq}_{m_{00}}\right)\cdot F_{1,{m_{00}}}\left(x^{\ell}\right)  =\frac{-t^{\ell}\cdot \left(\eta^{\ell,\geq}_{m_{00}}+\eta^{\ell,\leq}_{m_{00}}\right)}{F_{2,{m_{00}}}{\left(x^{\ell}\right)}} \le 0.
\]
By taking the limit here, we see that $\bar \varrho_{2,{m_{00}}} \leq 0$. Overall, we have shown that (\ref{eq:tstat-1})-(\ref{eq:tstat-3}) from Definition \ref{def:t-stat} hold at $\bar x$, hence, it is a T-stationary point.
\qed

\section{Application to SCNO}

We consider the sparsity constrained nonlinear optimization:
\[
\label{eq:SCNO}
\mbox{SCNO:}\quad \min_{x \in \R^n}\,\, f(x)\quad \mbox{s.\,t.} \quad \NNorm{x}{0} \le s,
\]
where the so-called $\ell_0$ "norm" counts non-zero entries of $x$:
\[
\NNorm{x}{0}= \AV{\left\{i\in \{1,\ldots,n\} \,\left\vert\, x_i \ne 0 \right.\right\}},
\]
the objective function $f \in C^2\left(\R^n,\R\right)$ is twice continuously differentiable, and $s \in\{0,1, \ldots,n-1\}$ is an integer.
For an SCNO feasible point $x$ we define the following complementary index sets:
\[
I_0(x) = \left\{i\in \{1,\ldots,n\} \,\left\vert\, x_i = 0 \right.\right\}, \quad
  I_1(x) = \left\{i\in \{1,\ldots,n\} \,\left\vert\, x_i \ne 0 \right.\right\}.
\]

Let us briefly recall the definition of M-stationarity for SCNO. 

\begin{definition}[M-stationary point, \cite{burdakov:2016}]
\label{def:m-stat}
A feasible point $\bar x$ is called M-stationary for SCNO if 
\[
\frac{\partial f}{\partial x_i} \left( \bar x\right) = 0  \mbox{ for all } i \in I_1\left(\bar x\right).
\]
\end{definition}
The concept of M-stationarity gains a lot of attention in the context of SCNO. 
In \cite{burdakov:2016}, the authors proved that M-stationarity is a necessary optimality condition for SCNO. Moreover, a relaxation method constructed there converges to an M-stationary point. In \cite{bucher:2018}, a second-order condition was stated, under which the local uniqueness of M-stationary points can be guaranteed.
In \cite{laemmel:2019},  the  concept  of  M-stationarity was used in  order  to  describe  the  global  structure  of  SCNO. By introducing nondegeneracy for M-stationary points, appropriate versions of deformation and cell-attachment theorems were shown to hold for SCNO. 

Our goal is now to examine the relaxation approach to SCNO as presented in \cite{burdakov:2016}. There, the following observation has been made: $\bar x$ solves SCNO if and only if there exists $\bar y$ such that $\left(\bar x, \bar y\right)$ solves the mixed-integer program:
\begin{equation}
   \label{eq:mix-int}
   \min_{x,y} \,\, f(x) \quad \mbox{s.\,t.} \quad 
   \sum_{i=1}^{n} y_i \geq n - s,  \quad 
      y_i \in \{0,1\},  \quad 
      x_iy_i =0,  \quad  i=1, \ldots, n.
\end{equation}
Using the standard relaxation of the binary constraints $y_i\in\{0,1\}$, the authors arrive at the  continuous optimization problem:
\begin{equation}
   \label{eq:relax}
   \min_{x,y} \,\,f(x) \quad \mbox{s.\,t.} \quad 
   \sum_{i=1}^{n} y_i \geq n - s,  \quad 
      y_i \in [0,1],  \quad 
      x_iy_i =0,  \quad  i=1, \ldots, n.
\end{equation}
 As pointed out in \cite{burdakov:2016}, SCNO and the optimization problem (\ref{eq:relax}) are closely related: $\bar x$ solves SCNO if and only if there exists a vector $\bar y$ such that $\left(\bar x, \bar y\right)$ solves (\ref{eq:relax}). 
 Obviously, the relaxation (\ref{eq:relax}) falls into the scope of MPOC due to the appearance of orthogonality type constraints. Hence, the notion of T-stationarity from Definition \ref{def:t-stat} can be applied to (\ref{eq:relax}). To do so, we use the following index sets given a feasible point $(x,y)$ of (\ref{eq:relax}):
\[
 \begin{array}{lcl}
\displaystyle  I_{00}\left( x, y\right) &=& \displaystyle  \left\{ i \in \{1,\ldots, n\} \,\left|\,  x_i = 0,  y_i =0 \right.\right\},
 \\
\displaystyle  I_{01}\left( x,  y\right) &=& \displaystyle  \left\{ i \in \{1,\ldots, n\} \,\left|\,  x_i = 0,  y_i =1 \right.\right\},
\\
\displaystyle  I_{0\not=}\left( x,  y\right) &=& \displaystyle  \left\{ i \in \{1,\ldots, n\} \,\left|\,  x_i = 0, y_i \not =0 \right.\right\}.
 \end{array}{}
\]
 
 \begin{definition}[T-stationarity for relaxation]
 A feasible point $(\bar x,\bar y)$ of (\ref{eq:relax}) is T-stationary if there exist multipliers:
\[
   \bar \mu, \bar \mu_i, i \in I_{01}\left(\bar x, \bar y\right), \bar \sigma_{1,i}, i \in I_{0\not=}\left(\bar x, \bar y\right), \bar \sigma_{2,i}, i \in I_1(\bar x), \bar \varrho_{1,i}, \bar \varrho_{2,i}, i \in I_{00}\left(\bar x, \bar y\right),
\]
such that the following conditions hold:
\begin{equation}
   \label{eq:tsrel-1} 
   \begin{array}{rcl}
   \left( \begin{array}{c}
        Df(\bar x)   \\
        0 
   \end{array}\right) &= & \displaystyle
   \bar \mu \left( \begin{array}{c}
        0  \\
        e 
   \end{array}\right) -
   \sum_{i \in I_{01}\left(\bar x, \bar y\right)} \bar \mu_i \left( \begin{array}{c}
        0  \\
        e_i 
   \end{array}\right) \\ && \displaystyle +
   \sum_{i \in I_{0\not=}\left(\bar x, \bar y\right) } \bar  \sigma_{1,i} \left( \begin{array}{c}
        e_i  \\
        0 
   \end{array}\right) + \sum_{i \in I_1(\bar x)} \bar \sigma_{2,i} \left( \begin{array}{c}
        0  \\
        e_i 
   \end{array}\right) \\ && \displaystyle +
   \sum_{i \in I_{00}\left(\bar x, \bar y\right)} \left(
   \bar \varrho_{1,i} \left( \begin{array}{c}
        e_i  \\
        0 
   \end{array}\right)  + \bar \varrho_{2,i} \left( \begin{array}{c}
        0  \\
        e_i 
   \end{array}\right) \right),
   \end{array}
\end{equation}
\begin{equation}
   \label{eq:tsrel-2}
   \bar \mu \left(\sum_{i \in I_0(\bar x)} \bar y_i - (n-s)\right) = 0 \mbox{ and } \bar \mu \geq 0, \quad \bar \mu_i \geq 0 \mbox{ for all } i \in I_{01}\left(\bar x, \bar y\right),
\end{equation}
\begin{equation}
   \label{eq:tsrel-3}
   \bar \varrho_{1,i}=0 \mbox{ or }\bar \varrho_{2,i}\le 0 \mbox{ for all } i \in I_{00}\left(\bar x, \bar y\right).
\end{equation}
 \end{definition}
 
 It turns out that T-stationary points of (\ref{eq:relax}) correspond to M-stationary points of SCNO.

\begin{proposition}[T- and M-stationarity]
\label{prop:tm-stat}
A feasible point $(\bar x, \bar y)$ is T-stationary for the relaxation (\ref{eq:relax}) if and only if the point $\bar x$ is M-stationary for SCNO. 
\end{proposition}

\proof Let $(\bar x, \bar y)$ be a T-stationary point of (\ref{eq:relax}). Since it is feasible for (\ref{eq:relax}), we have $\NNorm{\bar x}{0} \le s$. Otherwise, $\left|I_1(\bar x)\right| > s$, hence, $\bar y_i=0$ for $i \in I_1(\bar x)$. Together with $\bar y_i \leq 1$, $i \in I_0(\bar x)$ we then have:
\[
\sum_{i=1}^{n} \bar y_i = \sum_{i \in I_0(\bar x)} \bar y_i \leq \left|I_0(\bar x)\right| = n - \left|I_1(\bar x)\right| < n - s,
\]
a contradiction to the feasibility of $(\bar x,\bar y)$ for (\ref{eq:relax}).
Further, (\ref{eq:tsrel-1}) provides in particular:
\[  
    Df(\bar x) = \sum_{i \in I_{0\not=}\left(\bar x, \bar y\right) } \bar  \sigma_{1,i} e_i + \sum_{i \in I_{00}\left(\bar x, \bar y\right)} \bar \varrho_{1,i} e_i.
\]
Since $I_{0\not=}\left(\bar x, \bar y\right) \subset I_0(\bar x)$ and $I_{00}\left(\bar x, \bar y\right) \subset I_0(\bar x)$, we conclude that $\bar x$ is M-stationary for SCNO.

Vice versa, let $\left(\bar x, \bar y\right)$ be a feasible point of (\ref{eq:relax}) with $\bar x$ being M-stationary for SCNO. The multipliers can be chosen as follows in order to fulfil (\ref{eq:tsrel-1})-(\ref{eq:tsrel-3}):
\[
    \bar \mu=0, \quad \bar \mu_i = 0, i \in I_{01}\left(\bar x, \bar y\right), \quad \bar \sigma_{1,i} = \frac{\partial f}{\partial x_i}(\bar x), i \in I_{0\not=}\left(\bar x, \bar y\right), \quad \bar  \sigma_{2,i} = 0, i \in I_1(\bar x), 
\]
\[
    \bar \varrho_{1,i} = \frac{\partial f}{\partial x_i}(\bar x) \mbox{ and }
    \bar \varrho_{2,i} = 0, i \in I_{00}\left(\bar x, \bar y\right).
\]
Obviously, this choice of multipliers satisfies (\ref{eq:tsrel-2}) and (\ref{eq:tsrel-3}). In view of $I_{0\not=}\left(\bar x, \bar y\right) \cup I_{00}\left(\bar x, \bar y\right) = I_0(\bar x)$, also (\ref{eq:tsrel-1}) is fulfilled. Hence, $\left(\bar x, \bar y\right)$ is a T-stationary point of (\ref{eq:relax}). \qed

Proposition \ref{prop:tm-stat} states that by applying the proposed MPOC theory to (\ref{eq:relax}) we justified the use of M-stationarity for SCNO once again.
Let us mention that in \cite{burdakov:2016} the theory of  nonlinear programming was applied to the relaxation (\ref{eq:relax}).  

\begin{remark}[S-stationarity]
In \cite{burdakov:2016}, the notion of S-stationarity has been introduced for the relaxation (\ref{eq:relax}). By using our notation, a feasible point $\left(\bar x, \bar y\right)$ of (\ref{eq:relax}) is called S-stationary if there exist multipliers $\bar \sigma_{1,i}$, $i \in I_{0\not=}\left(\bar x, \bar y\right)$, such that
\[
 Df(\bar x) = \sum_{i \in I_{0\not=}\left(\bar x, \bar y\right) } \bar \sigma_{1,i} e_i.
\]
This definition is justified by the fact that a feasible point $\left(\bar x, \bar y\right)$ satisfies the Karush-Kuhn-Tucker condition if and only if it is S-stationary for (\ref{eq:relax}), see \cite{burdakov:2016}. We argue that any S-stationary point $\left(\bar x, \bar y\right)$ of (\ref{eq:relax}) is T-stationary. To see this we set:
\[
    \bar \mu=0, \quad \bar \mu_i = 0, i \in I_{01}\left(\bar x, \bar y\right), \quad \bar \sigma_{1,i} = \frac{\partial f}{\partial x_i}(\bar x), i \in I_{0\not=}\left(\bar x, \bar y\right), \quad \bar  \sigma_{2,i} = 0, i \in I_1(\bar x), 
\]
\[
    \bar \varrho_{1,i} = 0 \mbox{ and }
    \bar \varrho_{2,i} = 0, i \in I_{00}\left(\bar x, \bar y\right).
\]
These multipliers ensure the validity of conditions (\ref{eq:tsrel-1})-(\ref{eq:tsrel-3}). The reverse implication is not true in general, since for T-stationarity the multipliers $\bar \varrho_{1,i}$, $i \in I_{00}\left(\bar x, \bar y\right)$ do not necessarily vanish.  \qed
\end{remark}

 Let us now focus on the properties of T-stationary points for (\ref{eq:relax}). It turns out that all of them are degenerate in the sense of Definition \ref{def:nondeg}.
 
\begin{theorem}[Degeneracy]
\label{thm:deg}
All T-stationary points of the relaxation (\ref{eq:relax}) are degenerate.
\end{theorem}

\proof Let a T-stationary point $\left(\bar x, \bar y\right)$ of (\ref{eq:relax}) be given, i.\,e. (\ref{eq:tsrel-1})-(\ref{eq:tsrel-3}) are fulfilled. We show that at least one of the nondegeneracy conditions ND1-ND4 from Definition \ref{def:nondeg} is violated at $\left(\bar x, \bar y\right)$. For that, let us consider the following cases:

Case 1: $\displaystyle \sum_{i \in I_0(\bar x)} \bar y_i = n-s$. \\
If $I_{01}\left(\bar x, \bar y\right) \cup I_{00}\left(\bar x, \bar y\right) = I_0(\bar x)$, then LICQ is violated, i.\,e.~ND1 fails. Otherwise, the strict complementarity in ND2 is violated since $\bar \mu=0$, cf. (\ref{eq:tsrel-1}). 

Case 2: $\displaystyle \sum_{i \in I_0(\bar x)} \bar y_i > n-s$. \\
From (\ref{eq:tsrel-1}) we see that $\bar \mu_i=0$ for all $i \in I_{01}\left(\bar x, \bar y\right)$. This means that ND2 is violated unless $I_{01}\left(\bar x, \bar y\right) = \emptyset$.
Analogously, $\bar \varrho_{2,i}=0$ for all $i \in I_{00}\left(\bar x, \bar y\right)$. Again, in order to prevent the failure of ND3, $I_{00}\left(\bar x, \bar y\right) = \emptyset$ must hold. Overall, the tangent space from ND4 reads then:
\[
    \begin{array}{ll}
      T_{\left(\bar x,\bar y\right)}M(\bar x,\bar y) = \{ \left(\xi_x,\xi_y\right) \in \R^{2n} \,\,| &
        \left(\begin{array}{cc}
            e_i & 0  
        \end{array}\right) \left(\begin{array}{c}
             \xi_x  \\
             \xi_y 
        \end{array}\right)=0, i \in I_{0\not=}\left(\bar x, \bar y\right), \\ \\ &
     \left(\begin{array}{cc}
            0 & e_i  
        \end{array}\right) \left(\begin{array}{c}
             \xi_x  \\
             \xi_y 
        \end{array}\right)=0, i \in I_{1}\left(\bar x\right)
         \}.
     \end{array}
\]
The corresponding Lagrange function is
\[
   L(x,y) = f(x) - \sum_{i \in I_{0\not=}\left(\bar x, \bar y\right) } \bar  \sigma_{1,i} x_i - \sum_{i \in I_1(\bar x)} \bar \sigma_{2,i} y_i.
\]
Since $\left(\begin{array}{c}
            0 \\ e_i  
        \end{array}\right) \in T_{\left(\bar x,\bar y\right)}M(\bar x,\bar y)$ for all $i \in I_{0}\left(\bar x\right)$, its restricted Hessian becomes singular:
\[
  \left(\begin{array}{cc}
            0 & e_i  
        \end{array}\right) D^2 L\left(\bar x,\bar y\right) \left(\begin{array}{c}
            0 \\ e_i  
        \end{array}\right)= \left(\begin{array}{cc}
            0 & e_i  
        \end{array}\right)\left(\begin{array}{cc}
            D^2f\left(\bar x\right) & 0 \\ 
            0& 0  
        \end{array}\right) \left(\begin{array}{c}
            0 \\ e_i  
        \end{array}\right)=0.
\]
Thus, ND4 is violated. \qed

We conclude that the relaxation (\ref{eq:relax}) -- viewed as MPOC -- introduces degeneracy with regard to its T-stationary points. This is mainly due to the use of auxiliary $y$-variables needed to rewrite the sparsity constraint. Note that these degeneracies are intrinsic for (\ref{eq:relax}), i.\,e.~they appear even if we deal with originally nondegenerate M-stationary points for SCNO, see \cite{laemmel:2019} for the latter notion.

\bibliographystyle{apalike}
\bibliography{lit.bib}
\end{document}